\documentclass[reqno]{amsart}

\usepackage{enumerate}
\usepackage{ifpdf}
\usepackage{amsmath}
\usepackage{amsfonts}
\usepackage{amssymb}
\usepackage{amsthm}
\usepackage[hyperfootnotes=false]{hyperref}
\usepackage{setspace}
\usepackage{amsrefs}

\newcommand{\ignore}[1]{}

\renewcommand{\Im}{\operatorname{Im}}

\newcommand{\codim}{\operatorname{codim}}
\newcommand{\sing}{\operatorname{sing}}

\newcommand{\abs}[1]{\left\lvert {#1} \right\rvert}

\newcommand{\C}{{\mathbb{C}}}
\newcommand{\R}{{\mathbb{R}}}

\newcommand{\D}{{\mathbb{D}}}

\newcommand{\bB}{{\mathbb{B}}}

\newcommand{\bP}{{\mathbb{P}}}

\newcommand{\sA}{{\mathcal{A}}}

\newcommand{\sF}{{\mathcal{F}}}

\newcommand{\sO}{{\mathcal{O}}}

\newcommand{\sS}{{\mathcal{S}}}

\newcommand{\sZ}{{\mathcal{Z}}}

\newcommand{\rank}{\operatorname{rank}}

\newtheorem{thm}{Theorem}[section]

\newtheorem{prop}[thm]{Proposition}

\newtheorem{cor}[thm]{Corollary}

\newtheorem{lemma}[thm]{Lemma}

\theoremstyle{definition}
\newtheorem{defn}[thm]{Definition}
\newtheorem{example}[thm]{Example}

\theoremstyle{remark}
\newtheorem{remark}[thm]{Remark}

\author{Ji\v{r}\'i Lebl}
\thanks{The author was in part supported by NSF grant DMS 0900885.}
\address{Department of Mathematics, University of Illinois
at Urbana-Champaign, 
Urbana, IL 61801, USA}
\curraddr{Department of Mathematics, University of California
at San Diego, La Jolla, CA 92093-0112, USA}
\email{jlebl@math.uiuc.edu}
\date{February 27 2012}

\ifpdf
\hypersetup{
  pdftitle={Algebraic Levi-flat hypervarieties in complex projective space},
  pdfauthor={Jiri Lebl}
}
\fi

\title%
{Algebraic Levi-flat hypervarieties in complex projective space}

\begin{document}

%\doublespace

\begin{abstract}
We study singular real-analytic Levi-flat hypersurfaces in complex
projective space.
We define the rank of an algebraic Levi-flat hypersurface and study the
connections between rank, degree, and the type and size of the singularity.
In particular, we study degenerate singularities of algebraic
Levi-flat hypersurfaces.
We then give necessary and sufficient conditions for a Levi-flat
hypersurface to
be a pullback of a real-analytic curve in $\C$ via a meromorphic function.
Among other examples,
we construct a nonalgebraic semianalytic Levi-flat hypersurface with
compact leaves that is a perturbation of an algebraic Levi-flat variety.
\end{abstract}

\maketitle

%\enlargethispage{\baselineskip}

%%%%%%%%%%%%%%%%%%%%%%%%%%%%%%%%%%%%%%%%%%%%%%%%%%%%%%%%%%%%%%%%%%%%%%%%%%%%%

\section{Introduction} \label{section:intro}

The purpose of this paper is to organize some basic results on singular
Levi-flat
hypersurfaces in complex projective space.
First, we define and study several invariants of algebraic Levi-flat
hypersurfaces.  Second, we give
necessary and sufficient conditions for a real-analytic Levi-flat
hypersurface to be defined by a global
meromorphic function
and therefore algebraic.
Along the way we give several examples to
illustrate the phenomena encountered.

A real smooth hypersurface in a complex manifold
is said to be Levi-flat if it is pseudoconvex from
both sides.  If the hypersurface is real-analytic and nonsingular, then it is
classical that in suitable local coordinates, it can be represented by $\Im z_1
= 0$.
Therefore, there are no local holomorphic invariants.  The
situation is different if we allow the hypersurface to have singularities.
Local questions about singular Levi-flat
hypersurfaces have been previously studied
by Bedford~\cite{Bedford:flat}, Burns and Gong~\cite{burnsgong:flat},
and the author~\cites{Lebl:lfnm, Lebl:thesis}.
See also the books~\cites{BER:book, Boggess:CR, DAngelo:CR}
for the basic language and
background.

Let $\bP^n$ be the $n$-dimensional complex projective space.
Lins Neto proved \cite{linsneto:note} that no nonsingular real-analytic
Levi-flat
hypersurfaces exist in $\bP^n$, $n \geq 3$.
There have since been much work on generalizing this result further.
A different approach for the real-analytic case was taken by
Ni and Wolfson~\cite{NiWolfson:nolf}.
Siu~\cite{Siu:nolf3}, Cao and Shaw~\cite{CaoShaw:nolf3}, and
most recently Iordan and Matthey~\cite{IordanMatthey} improved the
regularity requirement.  The $n=2$ case was studied by Siu~\cite{Siu:nolf2}
and recently in the real-analytic setting by Ivashkovich~\cite{Ivashkovich}.

Singular real-analytic Levi-flat hypersurfaces, however, are a
different story and many such hypersurfaces exist.
Instead of hypersurface,
we will use the term \emph{hypervariety} for a codimension one subvariety
to emphasize the possibility of singularities, and to emphasize it is
a closed subvariety.
Also, unless specifically stated, subvarieties are analytic,
not necessarily algebraic.
Furthermore, unless specifically stated a subvariety is real-analytic.
Let $H \subset U \subset \C^k$ be a real-analytic hypervariety,
i.e.\@ a closed real subvariety of an open set $U$ of real codimension one.
Let $H^*$ be the set of points near which $H$ is a smooth hypersurface.
The hypervariety $H$
is said to be \emph{Levi-flat}, if it is Levi-flat at all points of $H^*$.
$H^*$ is foliated by complex hypersurfaces, and this foliation
is called the
\emph{Levi foliation}.
Any (real) algebraic Levi-flat hypervariety
in $\C^n$ can be extended to a Levi-flat hypervariety in~$\bP^n$.

One method to obtain algebraic Levi-flat hypervarieties in
$\bP^n$ is to take a rational function
$R \colon \bP^n \to \C$, and a real-algebraic one-dimensional subset $S \subset
\C$, and to consider the set $H = \overline{R^{-1}(S)}$.
Then $H$ is
an algebraic Levi-flat hypervariety. 
If $H$ is algebraic, then 
all the leaves of the Levi foliation must be compact.  A leaf is
said to be compact if the closure is a subvariety in $\bP^n$ of the same
dimension.
If a hypervariety is to be given by 
$H = \overline{R^{-1}(S)}$, then at least locally it must be given
by $\overline{F^{-1}(T)}$ for some local meromorphic function $F$
and some real-analytic set $T \subset \C$.  We will relax this condition
and suppose that $F$ is constant along the leaves of $H$.
It turns out that
these conditions are in fact 
sufficient.

\begin{thm} \label{mainthmalt}
Let $H \subset \bP^n$, $n \geq 2$, be an irreducible
Levi-flat hypervariety with infinitely many compact leaves.
Assume that for each $p \in \overline{H^*}$, there
exists a neighborhood $U$ of $p$ and a meromorphic function $F$ defined
on $U$ such that $F$ is constant along leaves of $H^*$.

Then, there exists a global rational function $R \colon \bP^n \to \C$
and
a real-algebraic one-dimensional
subset $S \subset \C$ such that $H \subset
\overline{R^{-1}(S)}$.
In particular, $H$ is semialgebraic; it is contained in an
algebraic Levi-flat
hypervariety.
\end{thm}

To prove the theorem we must find two objects.  We must find an algebraic
set $S \subset \C$ and the function $R$.  We find
a foliation of $\bP^n$ extending the Levi-foliation of $H$ by using a result of Lins Neto.
To find $R$
we apply
a
result of Darboux and generalized by Jouanolou, which
says
that a foliation of $\bP^n$
with infinitely many compact leaves has a rational first integral.
Next,
we find the $S \subset \C$ by proving Lemma~\ref{lfmeroalg},
which says that the image of $H$ under $R$ must essentially be
our algebraic curve $S$.

We really need to only study semianalytic sets.  A set is semianalytic if
it is locally constructed from
real-analytic sets by finite union, finite intersection, and complement.
For a hypervariety $H$, the set
$\overline{H^*}$ is semianalytic.
We will state a version of
Theorem~\ref{mainthmalt} for semianalytic hypersurfaces,
see Theorem~\ref{samainthm}.

The hypothesis of compact leaves seems necessary.  Example~\ref{ex:nonalg} is
a perturbation of an algebraic Levi-flat hypervariety of $\bP^2$,
which is again
Levi-flat, closed, semianalytic (thus
contained near each point in a real-analytic
subvariety), but not algebraic.  The leaves of this hypersurface are complex
hyperplanes, but do not extend to a foliation of $\bP^2$.
It also seems likely that a closure of a noncompact leaf of a 
foliation of $\bP^2$ could be semianalytic, though no such example is known to the
author.

Not all algebraic Levi-flat hypervarieties arise in the above way.  One
particular feature of hypervarieties defined using rational functions is
the existence of a degenerate singularity (in the sense of Segre varieties)
in dimension 2 or higher.
There are, however, algebraic Levi-flat hypervarieties
that do not have a degenerate
singularity as we will show in Example~\ref{ex:nodegen}.

To study algebraic Levi-flat hypervarieties we define their rank.
It is the rank of the Hermitian form of the defining bihomogeneous
polynomial.  Equivalently, the rank is the minimum number of holomorphic
polynomials needed to write the defining polynomial as a difference of squared
norms.
This definition of rank was used by the author together with D'Angelo to study
a seemingly unrelated problem in~\cite{DL:families}.  See
also the book by D'Angelo~\cite{DAngelo:CR} for further applications of this
circle of ideas.  For example, it is useful to write a defining equation
of a hypersurface as a squared norm to characterize the complex varieties
contained in the hypersurface.  We will also study a local
analytic version of the rank.

A simple argument shows that the dimension of the singular set of an algebraic
Levi-flat hypervariety in $\bP^n$ must be at least $2n-4$.  If the hypervariety
has no nondegenerate singularities, the dimension of the singular set must be of
maximal possible dimension, and the rank of the defining equation must also be
large compared to the dimension.  We will write $H_{sing}$ for the singular set
of $H$.  By the singular set
we mean the set of points near which $H$ is not a smooth
submanifold.  $H_{sing}$ is not in general equal to the complement of $H^*$ as
defined above, and is only a semialgebraic (or semianalytic if $H$ is not
algebraic) set.

It is standard
to also consider the algebraic singular set of an algebraic hypervariety.
The algebraic singular set is the set of points
where the defining polynomial has vanishing gradient.  The
algebraic singular set always contains the analytic singular set ($H_{sing}$ as
defined above), and the containment can be proper.
A classical example is $y^3+2x^2y-x^4 = 0$, which has no analytic singularities
but has an algebraic singularity at the origin, see~\cite{BCR:realalg}
for more on these issues.  Therefore, we also have the concept
of an algebraic degenerate singularity.  Any degenerate singularity is
an algebraic degenerate singularity, but not necessarily vice versa.
We summarize what we can say
about degenerate singularities in the following theorem.

\begin{thm} \label{degenthm}
Let $H \subset \bP^n$, $n \geq 2$,
be an algebraic Levi-flat hypervariety of rank~$r$.
\begin{enumerate}[(i)]
\item
If $r \leq n$ then there exists a complex subvariety $S \subset H$
of dimension at least $n-r$ such that every point in $S$
is an algebraic degenerate singularity of $H$.
\item
If $\dim H_{sing} < 2n-2$ then $H$ has a degenerate singularity.
\item
If $\dim H_{sing} = 2n-4$ then there is a complex subvariety $S \subset H$
of dimension $n-2$ such that every point in $S$
is a degenerate singularity of $H$.
\end{enumerate}
\end{thm}

The structure of this paper is as follows.  In \S~\ref{basicpropsec} we give
some standard basic results about real subvarieties of complex projective space
and Levi-flat hypervarieties in particular.
In \S~\ref{ranksec} we introduce and discuss
the rank of the hypersurface.
In
\S~\ref{degensingsec} we will prove Theorem~\ref{degenthm} and study the set of
degenerate singularities of an algebraic Levi-flat hypervariety.
In \S~\ref{algmerosec} we
study Levi-flat hypervarieties defined by meromorphic
functions.
In \S~\ref{lfandfolisec} we study holomorphic foliations induced by
Levi-flat hypervarieties and prove two alternate versions of
Theorem~\ref{mainthmalt}.
In \S~\ref{extfolsec} we prove that foliations extend from Levi-flat
hypervarieties even without compact leaves.
And finally in \S~\ref{compactleavessec} we study nonalgebraic Levi-flat
hypervarieties and semianalytic sets with compact leaves.

The author would like to thank 
Prof.\@ Xianghong Gong for suggesting to study singular Levi-flat hypersurfaces
in projective space.  The author would also like to thank the referee of an
earlier version of this paper for suggesting a simplification of the proof
of Theorem~\ref{mainthmalt}, and the referee of the current version for
suggestions on improving the organization of the paper.

%%%%%%%%%%%%%%%%%%%%%%%%%%%%%%%%%%%%%%%%%%%%%%%%%%%%%%%%%%%%%%%%%%%%%%%%%%%%%

\section{Basic properties} \label{basicpropsec}

Let $\sigma \colon \C^{n+1} \setminus\{0\} \to \bP^n$ be the natural projection.
Suppose $X$ is a real-analytic subvariety of $\bP^n$.  Define the
set $\tau(X)$ to be the set of points $z \in \C^{n+1}$ such that
$\sigma(z) \in X$ or $z = 0$.
A real-analytic subvariety $X \subset \bP^n$ is said to be \emph{algebraic}
if $X = \sigma(V)$ for some real-algebraic complex cone $V$ in $\C^{n+1}$.
A set $S$ is a complex cone when $p \in S$ implies $\lambda p \in S$ for
all $\lambda \in \C$.
We will
say that an algebraic subvariety $X \subset \bP^n$
is of degree $d$, if $d$ is the smallest
integer such that you need real polynomials of degree at most $d$ to define
$\tau(X)$.
We first establish some standard and easy to see properties of real-analytic
subvarieties of $\bP^n$.  By a bihomogeneous polynomial we mean a polynomial
that is separately homogeneous in $z$ and $\bar{z}$.  That is,
$P(t z, \overline{sz}) = t^{d/2} s^{d/2} P(z,\bar{z})$.  Thus, by
a degree $d$ bihomogeneous polynomial we mean of bi-degree $(d/2,d/2)$.

\begin{prop} \label{basicfacts}
Suppose
$X \subset \bP^n$ is a real-analytic subvariety.
\begin{enumerate}[(i)]
\item \label{basicra1}
$\tau(X) \setminus \{0\}$ is a real-analytic subvariety of
$\C^{n+1} \setminus \{0\}$.
\item \label{basicsubanal}
$\tau(X)$ is subanalytic.
\item \label{basicalg}
$X$ is 
algebraic if and only if
$\tau(X)$ is a real-analytic subvariety.
\item \label{basicbihom}
If $X$ is an irreducible algebraic hypervariety of
degree $d$, then 
$\tau(X)$ is defined by the vanishing of a single real valued
bihomogeneous polynomial of degree $d$ (bi-degree $(d/2,d/2)$).
\end{enumerate}
\end{prop}

\begin{proof}
To see \eqref{basicra1}, take homogeneous coordinates $[z_1:\cdots:z_{n+1}]$.
Fix e.g.\@ $z_1=1$ and find a set of defining functions $\rho_j$ for $X$
in some open set in the affine coordinates $z_2,\ldots,z_{n+1}$.  Let
$\tilde{\rho}_j(z_1,z_2,\cdots,z_{n+1}) =
\rho_j(z_2/z_1,\ldots,z_{n+1}/z_1)$ to be our defining equation
in some open subset of in $\C^{n+1} \setminus \{ z_1 = 0 \}$.

To see \eqref{basicsubanal} let again
$[z_1:\cdots:z_{n+1}]$ be the homogeneous coordinates, and let us
work in the chart where $z_1 \not= 0$.
Let $\tilde{X}$ be the subvariety in this chart.
Take the semianalytic set
$(\tilde{X} \cap \bB_n) \times \D$, where $\D \subset \C$ is the 
unit disc and $\bB_n \subset \C^n$ is the unit ball.  Define the
function $\varphi \colon \C^n \times \C \to \C^{n+1}$ by
$\varphi(w,\xi) = (\xi,\xi w)$.  This map takes 
$(\tilde{X} \cap \bB_n) \times \D$ to a subanalytic set $Y \subset \tau(X)$.
Furthermore, as $X$ is compact, then there are finitely many such 
charts and sets $Y_j$.   The germ of $\cup_j Y_j$ at the origin
agrees with the germ of $\tau(X)$ at the origin, which is what we needed to
prove.

One direction of \eqref{basicalg} is clear, the other is the same as in the
holomorphic case.  Let $\rho$ be a real-analytic function that is
zero on $\tau(X)$ near the origin.  Let $\rho = \sum_j \rho_j$
be the decomposition into homogeneous parts.  Take $t \in(-1,1)$
and note $\rho(t z) = \sum_j \rho_j (t z) =
 \sum_j t^j \rho_j (z)$.  For a fixed $z \in X$ we have
a power series that is identically zero.  Hence each $\rho_j$ must be zero
on $X$ and $X$ is therefore algebraic.

Finally, let us prove \eqref{basicbihom}.  Let $p$ be a real polynomial
vanishing on $\tau(X)$.  Write
\begin{equation}
p(z,\bar{z}) = \sum_{j,k} p_{jk} (z,\bar{z})
\end{equation}
where $p_{jk}$ is homogeneous of order $j$ in $z$ and
of order $k$ in $\bar{z}$.  Note that if $z \in \tau(X)$,
then $\lambda z \in \tau(X)$ for all $\lambda \in \C$.  Hence,
if $z \in \tau (X)$ then for all $\lambda$
\begin{equation}
0 = \sum_{jk} p_{jk} (\lambda z,\bar{\lambda} \bar{z})
= \sum_{jk} \lambda^j \bar{\lambda}^k p_{jk} (z,\bar{z}) .
\end{equation}
If we complexify $\lambda$ and $\bar{\lambda}$,
we get a polynomial in two variables
that is identically zero.  Therefore, $p_{jk} (z,\bar{z}) = 0$ for
all $j$ and $k$.

Since $\tau(X)$ is a real cone, it must be defined by single irreducible
real homogeneous polynomial of lowest degree.  If the defining
equation is not real homogeneous, then we can find a smaller degree
homogeneous polynomial vanishing on $\tau(X)$.  Call this polynomial
$p$ and write $p_{jk}$ as above.
Both the real and the imaginary parts of $p_{jk}$ must vanish on $X$,
hence we can write $p_{jk}(z,\bar{z}) = A(z,\bar{z}) p(z,\bar{z})$
for some (complex valued) polynomial $A$.  The degree of $p$
is equal to the degree of $p_{jk}$ and both are real homogeneous.
Plugging in $tz$ for $z$ and dividing by $t^k$ we notice that
$p_{jk}(z,\bar{z}) = A(tz,t\bar{z}) p(z,\bar{z})$ for all $t \in \R$, in
particular when $t = 0$.  Hence $p_{jk}$ is a constant times $p$ and of course
$p_{jk} = p$.  As $p$ was real valued we are done.  Notice also
that $j=k$.
\end{proof}

It is equally easy to see that any real polynomial in $\C^n$ can be 
made into a bihomogeneous polynomial in $\C^{n+1}$ and defines
a real subvariety of $\bP^n$.

We will be using the Segre variety to study Levi-flat hypervarieties.  Let $H
\subset U \subset \C^k$ be a real hypervariety defined by $\rho(z,\bar{z}) =
0$ for
some real-analytic function $\rho$ defined in $U$.
Let $\operatorname{conj}(U) = \{ z \mid \bar{z} \in U \}$.
Suppose that the power series for $\rho$ converges in
$U \times \operatorname{conj}(U)$ and hence we may complexify $\rho$.
We define
the \emph{Segre variety} $\Sigma_p$ as the set
\begin{equation}
\Sigma_p := \{ z \in U \mid \rho(z,\bar{p}) = 0 \}.
\end{equation}
The following is classical (see also for example~\cite{burnsgong:flat})
but we prove it here for completeness.

\begin{prop} \label{segrevarsmooth}
Suppose that $H \subset U \subset \C^k$
is a Levi-flat hypervariety ($U$ is small enough as above)
and $p \in H^*$.  Then one component $\Sigma_p'$
of $\Sigma_p$ agrees as a germ with
the leaf of the Levi foliation of $H$ through $p$.  The germ
of $\Sigma_p'$ is the unique germ of a complex hypervariety through $p$.
\end{prop}

\begin{proof}
Taking $U$ smaller can at most make $\Sigma_p$ smaller.  Thus we can make $U$
small enough such that we can make a local change of coordinates such that
$H$ is given near $p$ as $\{ \Im z_1 = 0 \}$ and $p$ is the origin.  Then
$\rho(z,\bar{z}) = a(z,\bar{z})(1/2i)(z_1-\bar{z}_1)$.  $\Sigma_0$ then
contains $\{ z_1 = 0 \}$.  Let $f(z)$ be another holomorphic function such
that $\{ f = 0 \} \subset H$ and $f(0) = 0$, then $f$ is real valued on
$H$ and in particular on $\{ z_1 = 0 \}$, hence $f = 0$ when $z_1 = 0$, and
uniqueness follows.
\end{proof}

In particular, if $H$ is a Levi-flat hypervariety
and for some $p \in H$ there are two distinct germs
of a complex hypervariety through $p$ contained in $H$, then $p$ must be
a singular point of $H$.  Similarly, if there is a germ of a 
complex analytic hypervariety contained in $H$ through $p$, singular at $p$,
then $H$ itself must be singular at $p$.
If $H$ is of a higher codimension
at $p$, then the above two statements are obvious.

We now focus on Levi-flat hypervarieties $H \subset \bP^n$.
We get the following corollary.  We say a leaf $L$ of the Levi foliation
of $H^*$ is compact if the closure $\bar{L}$ has the same dimension.
In this case, by Remmert and Stein $\bar{L}$ is a complex subvariety.  We will
generally abuse terminology and call $\bar{L}$ a leaf of $H$.

\begin{cor}
Let $H \subset \bP^n$ be an algebraic Levi-flat hypervariety.  Then
all leaves of $H$ are compact.
\end{cor}

\begin{proof}
Take a defining polynomial for $\tau(H)$ and look at the
Segre varieties.  Note that any leaf
is either contained in some Segre variety that is proper subset of
$\C^{n+1}$, in which case it is a compact leaf, or it is contained in the set
of points where the Segre variety is not a proper subset of $\C^{n+1}$ (what
we will call the degenerate singularities).  But that subset itself must be a
proper complex subvariety.
\end{proof}

We have the following
simple, and surely classical, observation.

\begin{prop} \label{coneisolthen2}
If $H \subset \bP^n$ is an algebraic Levi-flat hypervariety.  Then
\begin{equation} \label{singineq}
2n-4 \leq \dim H_{sing} \leq 2n-2 .
\end{equation}
\end{prop}

\begin{proof}
Any leaf of $\tau(H)$ must pass through the origin
since $\tau(H)$ is a complex cone.
Hence, any two leaves must meet on a complex subvariety of dimension $n-1$
in $\C^{n+1}$
and this set must lie in the singular set.
\end{proof}

The singularity can also be of larger dimension than $2n-4$.
Pick any singular algebraic
curve $S$ in $\bP^1$ and look at $\tau(S) \subset \C^2$.
The singular set is going to
be a finite union of complex lines through the origin.
Of course this argument also implies that if $n \geq 2$,
then $H$ must be singular.

The canonical
local example of a Levi-flat hypersurface in $\C^n$ is defined
by $\Im z_1 = 0$.  This hypersurface can of course be extended to all of
$\bP^n$.  If we bihomogenize this equation we will get a quadratic
complex cone in $\C^{n+1}$ given by
\begin{equation} \label{quadraticcone}
z_1 \bar{z}_2 - \bar{z}_1 z_2 = 0 .
\end{equation}

Burns and Gong~\cite{burnsgong:flat} have classified, up to local
biholomorphism, all germs of quadratic Levi-flat hypervarieties.
I.e.\@ up to
biholomorphism, there is only one quadratic complex cone that is a Levi-flat
hypervariety,
and that is given by \eqref{quadraticcone}.
It is not hard to show this fact directly 
using Proposition~\ref{basicfacts} and it is equivalent to
the following proposition.

\begin{prop} \label{quadconeprop}
Suppose that $H$ is a quadratic algebraic Levi-flat hypervariety in $\bP^n$.
Then $H$ is biholomorphically equivalent to a hypervariety
given by \eqref{quadraticcone}.
\end{prop}

\begin{proof}
Let $\rho$ be the defining
bihomogeneous polynomial of degree 2 for $\tau(H)$.
Since it is of degree 2 and
bihomogeneous it can be written as a Hermitian form, i.e.\@ 
$\rho(z,\bar{z}) = \bar{z}^t A z$, where $z = (z_1,\ldots,z_{n+1})^t$.
As $A$ is Hermitian, we can make a linear change of variables
(a biholomorphic transformation of $\bP^n$) such that $A$ is diagonal.
Thus we can assume
\begin{equation}
\rho(z,\bar{z}) = \sum_j \epsilon_j z_j\bar{z}_j .
\end{equation}
We can make further linear transformations to assume that
$\epsilon_j = -1, 0,$ or $1$.  Being Levi-flat is equivalent to the
Levi form vanishing at all smooth points.  This condition
is equivalent to the following differential equation
\begin{equation} \label{flateqnice}
\rank
\left[
\begin{matrix}
\rho & \rho_z \\
\rho_{\bar{z}} & \rho_{z\bar{z}}
\end{matrix}
\right]
\leq 2
\ \ \ \text{ on $\rho = 0$ } .
\end{equation}
All 3 by 3 subdeterminants of the matrix
must be zero, hence all but two
$\epsilon_j$ must be zero.  It is not hard to see that at least two
must be nonzero and of different sign, otherwise $H$ is not a hypersurface.
Thus we can assume that $\rho(z,\bar{z}) = z_1 \bar{z}_1 - z_2 \bar{z}_2$,
which is unitarily equivalent to~\eqref{quadraticcone}.
\end{proof}

Therefore,
there exist affine coordinates such that every quadratic Levi-flat
hypervariety of $\bP^n$ is given by $\Im z_1 = 0$ in those affine coordinates.
When $n \geq 2$, the hypervariety is singular and the singular set is
a complex subvariety of dimension $n-2$.

We end the section with an example which illustrates the subtlety of 
the geometry of the singular set.
Levi-flat hypervarieties
generally suffer from the same subtle issues as do real-analytic subvarieties
in general.

\begin{example} \label{cartanumbex}
A classical example is the subvariety given by $y^2+x^2-x^3 = 0$ in $\R^2$.
We get an irreducible (algebraically)
curve for which the origin is an isolated point.  We can think of
$\R^2$ as $\C$ using $x+iy = z$ and let $X$ be the subvariety extended
to $\bP^1$.
The equation then becomes
\begin{equation}
\bar{z}^3+3z\bar{z}^2+3z^2\bar{z}-8z\bar{z}+z^3 = 0 .
\end{equation}
We bihomogenize this equation to get a complex cone in $\C^2$.
That is, we define the hypervariety $H = \tau(X)$ by
\begin{equation}
w^3\bar{z}^3+3w^2\bar{w}z\bar{z}^2+3w\bar{w}^2z^2\bar{z}
-8w^2\bar{w}^2z\bar{z}+\bar{w}^3z^3
= 0 .
\end{equation}
The left hand side is irreducible as a polynomial,
but also analytically at the origin.
Suppose $f$ is a real-analytic function defined on a neighborhood
of the origin that vanishes on a nontrivial part of $H^*$. 
Write $f = \sum f_j$ where $f_j$ are real homogeneous of degree $j$.
Using the proof
of Proposition~\ref{basicfacts} we see that each $f_j$ vanishes on a nontrivial
part of $H^*$.  Thus each $f_j$ vanishes on all of $H$,
as $H$ was defined by an irreducible polynomial.  Hence, $f$ vanishes on
all of $H$.

$H$ is Levi-flat as it is a complex cone in $\C^2$.
Near all points of the set $\{ z = 0, w \not= 0 \}$, $H$ is
a complex line.  Therefore $H^*$ does not include the set
$\{ z = 0, w \not= 0 \}$.  This set is colloquially called the ``stick''
of the ``umbrella.''

Note also, that the one-dimensional part of $X$ is a real-analytic subvariety,
but it is only semialgebraic.  So we have an example of a real-analytic
hypervariety of $\bP^1$, that is semialgebraic and not algebraic.
\end{example}

The ``stick'' of the umbrella need not be complex analytic.
Brunella~\cite{Brunella:lf} gives the following example.
Let $z = x+iy$ and $w = s+it$.  Then the set given by
$t^2 = 4(y^2 + s)y^2$
is Levi-flat and the ``stick'' is the set $\{ t = s = 0, s \leq 0 \}$,
which is totally real in $\C^2$.

Finally, we will need the following lemma of
Burns and Gong (Lemma~2.2 in \cite{burnsgong:flat}) to see that for
an irreducible hypervariety $H$ we need only require it to be Levi-flat
at one point of $H^*$.

\begin{lemma}[Burns-Gong] \label{burnsgonglemma}
Let $H \subset \C^n$ be a real-analytic hypervariety, locally irreducible
at point $p \in H$.  Then there exists an open set $U \subset \C^n$ containing
$p$ such that $H \cap U$ is Levi-flat if and only if one of the components
of $H^* \cap U$ is Levi-flat.
\end{lemma}

The lemma is essentially proved by noticing that equation \eqref{flateqnice}
must hold everywhere on $H$ by complexification
of the irreducible defining function of $H$.
By noticing again that the property of being Levi-flat is equivalent to
the equation \eqref{flateqnice},
we have the following trivial classical
proposition which is useful in
application of the Burns-Gong lemma.

\begin{prop}
Let $M \subset \C^n$ be a connected
real-analytic submanifold of real dimension
$2n-1$.  $M$ is Levi-flat if and only if there exists
an open set $N \subset M$ such that $N$ is Levi-flat.
\end{prop}

%%%%%%%%%%%%%%%%%%%%%%%%%%%%%%%%%%%%%%%%%%%%%%%%%%%%%%%%%%%%%%%%%%%%%%%%%%%%%

\section{Rank}
\label{ranksec}

Let $H \subset \bP^n$ be an algebraic Levi-flat hypervariety.
Let $P$ be
the defining bihomogeneous polynomial for $\tau(H)$.  Using multi-index notation
we write
\begin{equation}
P(z,\bar{z}) = \sum_{\alpha \beta} c_{\alpha \beta} z^\alpha \bar{z}^{\beta} .
\end{equation}
Hence, if we order the multi-indices in some way and write the column vector
$\sZ = (z^{\alpha_1}, z^{\alpha_2}, \ldots, z^{\alpha_m})^t$, we can write
the matrix $C = \left[ c_{\alpha \beta} \right]_{\alpha\beta}$, and then
\begin{equation}
P(z,\bar{z}) = \overline{\sZ}^t C \sZ
\end{equation}
As $P$ is real valued, then $c_{\alpha\beta} = \overline{c_{\beta\alpha}}$,
hence $C$ is Hermitian.

\begin{defn}
Let $H \subset \bP^n$ be a real-algebraic hypervariety and $P$
the defining polynomial for $\tau(H)$.  We form the matrix $C$ and
define
\begin{align}
& \rank P := \rank C ,\\
& \rank H := \rank C .
\end{align}
\end{defn}

It is standard that if $C$ is of rank $r$,
then there exist $r$ column vectors $v_1, v_2, \ldots, v_r$ such that
\begin{equation}
C =
v_1 \overline{v_1}^t  + \cdots +
v_2 \overline{v_s}^t -
v_{s+1} \overline{v_{s+1}}^t - \cdots -
v_r \overline{v_r}^t .
\end{equation}
Taking $p_j(z) := \overline{v_j}^t \sZ$ we can see that
\begin{equation} \label{Pdiffsq}
P(z,\bar{z}) =
\abs{p_1(z)}^2 + \cdots +
\abs{p_s(z)}^2 -
\abs{p_{s+1}(z)}^2 - \cdots -
\abs{p_r(z)}^2 .
\end{equation}
The number $r$ is the minimum number of
holomorphic polynomials $p_j$ we will need for
such a decomposition.
Hence, the rank $r$ can be also defined as the minimal number of
holomorphic polynomials such that $P$ can be written as \eqref{Pdiffsq}.

As $H$ is a hypersurface, there must be at least some
positive and some negative eigenvalues of $C$.  Therefore, $\rank H \geq 2$.
On the other hand, we have the trivial estimate
$\rank H \leq \binom{d/2+ n}{n}$.

\begin{prop} \label{rank2leviflatprop}
If $\rank H = 2$, then $H$ is Levi-flat.
\end{prop}

\begin{proof}
$H$ is the set $\abs{p_1(z)}^2 - \abs{p_2(z)}^2 = 0$, and hence
a Levi-flat hypervariety defined by the meromorphic
function $p_1 / p_2$.
\end{proof}

\begin{example} \label{lfhighrank}
Of course there exist Levi-flat hypervarieties with higher rank.  For example,
Let $z = x+iy$.  The real curve $x^3-y^2 = 0$ in $\C$ can be extended to
$\bP^1$ (or $\bP^n$ by considering the equation in $\C^n$)
by bihomogenizing the equation (using the variable $w$) to get
the polynomial
\begin{equation}
w^3 \bar{z}^3
+ 3 z w^2 \bar{z}^2 \bar{w}
+ 2 w^3 \bar{z}^2 \bar{w}
+ 3 z^2 w \bar{z} \bar{w}^2
- 4 z w^2 \bar{z} \bar{w}^2
+ z^3 \bar{w}^3
+ 2 z^2 w \bar{w}^3 .
\end{equation}
If we let $\sZ = (z^3, z^2 w, z w^2, w^3)^t$, we get the
matrix
\begin{equation}
C = 
\begin{bmatrix}
0 & 0 & 0 & 1 \\
0 & 0 & 3 & 2 \\
0 & 3 & -4 & 0 \\
1 & 2 & 0 & 0
\end{bmatrix} .
\end{equation}
The rank is 4, there are 2 positive and 2 negative eigenvalues.  We can use the identity $a \bar{b} + \bar{a} b = 
\abs{a + b}^2 - \abs{a-b}^2$ to actually find a decomposition of the
polynomial
as follows
\begin{multline}
\abs{z^3+2z^2w + w^3}^2 - 
\abs{z^3+2z^2w - w^3}^2
+
\\
\abs{3z^2w-2zw^2 + zw^2}^2 - 
\abs{3z^2w-2zw^2 - zw^2}^2 .
\end{multline}
Hence, we have found an example where the rank of $H$ is equal to the
maximal rank possible and $H$ is still Levi-flat.
\end{example}

\begin{remark} \label{rankremark}
Note that Proposition~\ref{quadconeprop} also says that any quadratic Levi-flat
hypervariety of $\bP^n$ must have rank 2.  On the other hand, the generic
quadratic hypervariety in $\bP^n$ has rank $n+1$.
Hence, it is not always
possible to construct examples that are Levi-flat and have the maximal rank
$\binom{d/2+ n}{n}$.  However,
Example~\ref{ex:nodegen} is a Levi-flat hypervariety
of $\bP^2$ of degree 4 and has rank 6, which is the maximum possible.

The above example however gives a way to construct examples of arbitrarily high
rank.  Given any real-algebraic curve in $\C$ we can bihomogenize
the defining equation and get a Levi-flat
cone in $\C^{n+1}$, for any $n \geq 1$, and
thus get a Levi-flat hypervariety of $\bP^n$.  As we can choose an irreducible
curve of degree $\delta$ such that
the bihomogenized polynomial can have rank $\delta+1$, the maximal possible
rank for a curve in $\C$.
\end{remark}

\begin{prop} \label{rankinvar}
Let $H \subset \bP^n$ be a real-algebraic hypervariety.
Then $\rank H$ is invariant under automorphisms
of $\bP^n$.
\end{prop}

\begin{proof}
Let $L$ be an invertible linear
mapping and $P$ is given by \eqref{Pdiffsq}, then $L^{-1}(\tau(H))$
is given by
\begin{equation} \label{PdiffsqL}
P(Lz,\overline{Lz}) =
\abs{p_1(Lz)}^2 + \cdots +
\abs{p_s(Lz)}^2 -
\abs{p_{s+1}(Lz)}^2 - \cdots -
\abs{p_r(Lz)}^2 .
\end{equation}
Thus the rank cannot increase (and therefore cannot decrease) by composing
with a linear transformation $\C^{n+1}$ and hence an automorphism of $\bP^n$.
\end{proof}

It is also possible to work in some set of affine coordinates, rather than the
homogeneous coordinates.  The rank can be defined in generic
affine coordinates and
we will get the same number as we get in homogeneous coordinates.
This procedure
suggests that we might similarly define the rank (locally) for a nonalgebraic
hypervariety.  We get a genuinely different notion of rank, which we study in
the next section.

\section{Degenerate singularities}
\label{degensingsec}

Let $H$ be a Levi-flat hypervariety defined near $p \in \C^k$.
For each defining function $\rho$ of $H$ we find a neighborhood
$U$ of $p$ small enough such that $\rho$ complexifies
as in \S~\ref{basicpropsec} and we may define the Segre variety
$\Sigma_q$
for all $q \in U$.

\begin{defn}
We will say that $p \in H$
is a \emph{degenerate singularity} if the Segre variety
$\Sigma_p$ is open (of dimension $k$) for every local defining 
function.  In other words, $p$ is a degenerate singularity if
$z \mapsto \rho(z,\bar{p})$
is identically zero for $z$ near $p$ for
all local defining functions $\rho$ of $H$.

If $H$ is algebraic and $P(z,\bar{z})$ is the defining polynomial of $H$,
then an algebraic singular point $p \in H$ (a point where
the gradient of $P$ vanishes) is called
an \emph{algebraic degenerate singularity} of $H$
if $z \mapsto P(z,\bar{p})$ is identically zero.
\end{defn}

We fix a local defining function $\rho$
and a connected neighborhood $U$ as above.
By the reality of $\rho$ we note that if $q \in \Sigma_p$ then
$p \in \Sigma_q$.  Hence, if $p$ is a degenerate singularity,
then $p \in \Sigma_q$ for all $q \in U$.  On the other hand if
$q \in \Sigma_p$ for all $q \in H \cap U$, then $q \in \Sigma_p$
for all $q \in U$ and hence $p$ is degenerate.
Therefore by
Proposition~\ref{segrevarsmooth} we have the following result.

\begin{prop}
If $p$ is a singularity of a Levi-flat hypervariety $H$
and there are infinitely many distinct
germs of complex hypervarieties $(L,p) \subset (H,p)$, then
$p$ is a degenerate singularity.
\end{prop}

The hypervariety defined by 
\begin{equation} \label{quadconedegsing}
z_1 \bar{z}_2 - \bar{z}_1 z_2 = 0
\end{equation}
has a degenerate singularity at 0.  By the reasoning above, whenever
$H$ is a complex cone, the origin
is always a degenerate singularity.  

Let $P(z,\bar{z})$ be a defining polynomial for a real-algebraic 
hypervariety $H \subset \C^k$.  We call the Segre variety 
induced by $P$ the \emph{algebraic Segre variety}.
Then a singular point $p$ is an algebraic
degenerate singularity if the algebraic Segre variety $\Sigma_p = \C^k$.
If a point is a
degenerate singularity, then it must also be 
an algebraic degenerate singularity.

\begin{remark}
It may happen that the Segre variety for the local analytic defining function is
different from the one given by the defining polynomial,
as is illustrated by the following two examples.
Let $z_1 = x+iy$, $z_2 = s+it$.
The classical example
\begin{equation}
y^3+2x^2y-x^4 = 0
\end{equation}
is a Levi-flat
hypervariety with algebraic (but not analytic) singularity along the
set $\{ z_1 = 0 \}$.  The algebraic Segre variety is a triple plane at
the origin union another disjoint complex hyperplane, while the analytic
Segre variety is only a small piece of $\{ z_1 = 0 \}$.  To see a more
dramatic example consider the equation
\begin{equation}
s^3+2x^2s-x^4 = 0 .
\end{equation}
This equation defines a
real-algebraic hypervariety (not Levi-flat) with
an algebraic (but not analytic) singularity on the set $\{ x = s = 0 \}$.
The hypervariety is a smooth real-analytic submanifold and so the
Segre variety induced by the local analytic defining equation is a nonsingular
complex hypersurface.  However, the algebraic Segre variety is locally
a union of 3 smooth complex hypersurfaces at the origin.
\end{remark}

We have found in Proposition~\ref{quadconeprop}
that the only quadratic algebraic Levi-flat hypervariety
in $\bP^n$ is the quadratic cone given in homogeneous coordinates
by \eqref{quadconedegsing}.  When $n \geq 2$,
this hypervariety always has a degenerate singularity.  On the other hand
it is sufficient to consider degree 4 in $\bP^2$ to have
an example without a degenerate singularity.

\begin{example} \label{ex:nodegen}
To construct
a Levi-flat hypervariety of $\bP^2$
without any (algebraic or analytic) degenerate singularities,
we 
construct a Levi-flat hypervariety $H \subset \C^3$ that is a complex
cone and such that the origin is the only degenerate singularity.
We look at the equation
\begin{equation} \label{eq:planeeq}
z_1 + z_2 t + z_3 t^2 = 0 .
\end{equation}
We look the points $z = (z_1,z_2,z_3)$ where this
polynomial has a real solution.  We thus have the semialgebraic surface
\begin{equation}
\{ z \in \C^3 \mid z_1 + z_2 t + z_3 t^2 = 0 \text{ for some } t \in \R \}.
\end{equation}
By applying the quadratic formula and finding where the solution is real,
we obtain the following degree 4 real homogeneous polynomial defining a
Levi-flat hypervariety, which contains all the planes defined
in~\eqref{eq:planeeq}.
\begin{equation} \label{eq:planeeqfull}
z_1^2 \bar{z}_3^2
+ z_1 z_2 \bar{z}_2 \bar{z}_3
+ z_2^2 \bar{z}_1 \bar{z}_3
+ z_1 z_3 \bar{z}_2^2 - 2 z_1 z_3 \bar{z}_1 \bar{z}_3
+ z_2 z_3 \bar{z}_1 \bar{z}_2
+ z_3^2 \bar{z}_1^2
= 0 .
\end{equation}
The hypervariety is Levi-flat
by applying 
Lemma~\ref{burnsgonglemma} and noticing that
at least some subset of the hypervariety is
foliation
by the varieties defined by
$z_1 + z_2 t + z_3 t^2 = 0$ for some fixed real $t$.
The only point that lies in all
these varieties is the origin.  It is easy to see that even the entire
hypervariety has no degenerate singularities (except the origin)
by writing the defining equation as
\begin{equation}
z_1^2 \big( \bar{z}_3^2 \big)
+ z_1 z_2 \big( \bar{z}_2 \bar{z}_3 \big)
+ z_2^2 \big( \bar{z}_1 \bar{z}_3 \big)
+ z_1 z_3 \big ( \bar{z}_2^2 - 2 \bar{z}_1 \bar{z}_3 \big)
+ z_2 z_3 \big( \bar{z}_1 \bar{z}_2 \big)
+ z_3^2 \big( \bar{z}_1^2 \big)
= 0 .
\end{equation}
We think of $z$ and $\bar{z}$ as independent.  To find the (algebraic)
Segre variety corresponding to this equation we would set $\bar{z}$
to a constant.
So the expressions in parentheses are coefficients of a polynomial
in $z$.  The only place where they all vanish identically is 
when $\bar{z}_1 = \bar{z}_2 = \bar{z}_3 = 0$.
Thus $\sigma(H) \subset \bP^2$
has no algebraic degenerate singularities, and hence no analytic degenerate
singularities.  The variety $\sigma(H)$ is of rank 6 as is easily seen
from the defining equation~\eqref{eq:planeeqfull}.
\end{example}

The size of the algebraic degenerate singular set is related to the
rank of $H$, which is itself related to the degree.
The following lemma proves
the first part of
Theorem~\ref{degenthm} from the introduction.

\begin{lemma} \label{degenranklemma}
Let $H \subset \bP^n$, $n \geq 2$,
be an algebraic Levi-flat hypervariety of rank $r$.
If $r \leq n$ then there exists a complex subvariety $S \subset H$
of dimension at least $n-r$ such that every point in $S$ is an algebraic
degenerate singularity of $H$.
\end{lemma}

In particular, if $H$ is nondegenerate then $\rank H > n$.  The proof
of this lemma is essentially the following observation, which we state as a
proposition.  This result essentially gives us a method to find all algebraic
degenerate singularities of $H$.

\begin{prop} \label{degenrankprop}
Let $H$ be as above and $P$ the defining bihomogeneous polynomial,
and let $r$ be the rank.  Write
\begin{equation} \label{Pdiffsq2}
P(z,\bar{z}) =
\abs{p_1(z)}^2 + \cdots +
\abs{p_s(z)}^2 -
\abs{p_{s+1}(z)}^2 - \cdots -
\abs{p_r(z)}^2 .
\end{equation}
Then $w$ is an algebraic degenerate singularity of $\tau(H)$ if
and only if $p_j(w) = 0$ for all~$j = 1,\ldots,r$.
\end{prop}

Note that if $p_j(w) = 0$ for all $j$, then $w \in \tau(H)$,
and $z \mapsto P(z,\bar{w})$ is identically zero.  As the rank is $r$,
the $p_j$ are linearly independent.  The converse then follows.
We finish the proof of Lemma~\ref{degenranklemma}
by applying Proposition~\ref{degenrankprop}
and taking $S$ to be the subvariety defined by
$p_j = 0$ for all $j$.

It is not true that high rank guarantees lack of degeneracy.  Since any
real-algebraic curve in $\C$ extends to a Levi-flat hypervariety
in $\bP^n$ as in Remark~\ref{rankremark} we can get Levi-flat hypersurfaces
with arbitrarily high rank.  However, all hypervarieties obtained in this
way will have degenerate singularities.

By Proposition~\ref{coneisolthen2} we note that the singular set of
an algebraic Levi-flat hypervariety of $\bP^n$ has to be
at least of real dimension $2n-4$.
This fact follows because when two leaves of the foliation meet,
they must meet in
a set of complex dimension $n-2$ and this set must be contained in the
singular set of $H$.  It therefore easy to see that if the singular set is
only of dimension $2n-4$ then all the leaves must meet on the same set of
dimension $n-2$.  Hence we have the following lemma, which proves the
second part of Theorem~\ref{degenthm}.

\begin{lemma}
Let $H \subset \bP^n$, $n \geq 2$, be an algebraic Levi-flat hypervariety
such that $\dim H_{sing} = 2n-4$.  Then there exists a complex subvariety
$S$ of dimension $n-2$ such that every point in $S$ is a degenerate 
singularity of $H$.
\end{lemma}

Obviously if $H$ is not to have any degenerate singularity, then 
the singular set must be large.
The author essentially
proved in~\cite{Lebl:lfnm}
that if the singular set is a submanifold of dimension $2n-2$,
in the hypervariety, it is either complex or Levi-flat (i.e.\@ locally 
equivalent to $\C^{n-2} \times \R^2$).
The following lemma tells us that 
nondegeneracy must be
compensated by such a singular set and proves the final
part of Theorem~\ref{degenthm}.

\begin{lemma}
If $H \subset \bP^n$, $n \geq 2$, is an algebraic
Levi-flat hypervariety
without degenerate singularities, then the singular set must be of
real dimension $2n-2$.
\end{lemma}

\begin{proof}
We look at $\tau(H) \subset \C^{n+1}$.  We look at the leaves 
of the Levi foliation going through the origin.  Any two
such leaves
must meet on a set of complex dimension $n-1$, and this set must lie in the
singular set of $\tau(H)$.
As before, if all leaves met on the same set, then 
$\tau(H)$ would have a degenerate singularity away from the origin and
hence $H$ would have a degenerate singularity.
Thus suppose that the singular set of $\tau(H)$ is of dimension
$2n-1$.
Let us look at
a family of leaves $\{ L_t \}$ parametrized by a real parameter $t$ in some
small interval $(-\epsilon,\epsilon)$.  That is, find
\begin{equation}
P(z,t) := \sum_{\abs{\alpha} = d} a_\alpha (t) z^\alpha ,
\end{equation}
where $a_\alpha(t)$ are real-analytic functions in $t$, and
such that the sets $L_t = \{ z | P(z,t) = 0 \}$ are leaves of $\tau(H)$.
We can find such a $P$ by considering the coefficients of a polynomial
in $z$ as variables and then the set of polynomials whose
zero sets are contained in $\tau(H)$ is a semialgebraic set.

Take two such parameters and look
at the set $L_t \cap L_s$.
The sets 
$L_t \cap L_s$ have real dimension $2n-2$ (complex dimension $n-1$).
Fix $t$ and note
$\bigcap_{s \not= t} L_s \cap L_t = \emptyset$.
Hence, there must exist a submanifold $T_t \subset \tau(H)_{sing}$
of dimension $2n-1$ that is
foliated by $(n-1)$-dimensional complex submanifolds (the $L_t \cap L_s$).
We can pick a maximal such $T_t$ (not necessarily unique).

For each $t$ such a statement is true and as the
singular set is of dimension $2n-1$,
there is
some $t_0$ such that for infinitely $t$, the set
$T_{t_0} \cap T_{t}$ is nonempty
and of dimension $2n-1$.  But then infinitely many $L_t$
have the same nontrivial intersection with $L_{t_0}$, and hence
the hypersurface would
have a degenerate singularity.  We obtain a contradiction.
Consequently, the singular set of $\tau(H)$ must be of dimension $2n$,
and
the singular set of $H$ was real $(2n-2)$-dimensional.
\end{proof}

%%%%%%%%%%%%%%%%%%%%%%%%%%%%%%%%%%%%%%%%%%%%%%%%%%%%%%%%%%%%%%%%%%%%%%%%%%%%%

\section{Algebraic Levi-flat hypervarieties defined by meromorphic functions}
\label{algmerosec}

The following construction
gives a large supply of Levi-flat hypervarieties of $\bP^n$,
although it is not exhaustive.

\begin{prop}
Let $F$ be a meromorphic function on $\bP^n$.
Let $S \subset \C$ be a real-algebraic curve.
Then the set $\overline{F^{-1}(S)}$
is an
algebraic Levi-flat hypervariety.  However, not all
algebraic Levi-flat hypervarieties in $\bP^n$ are defined in this manner.
\end{prop}

\begin{proof}
It is standard that
any meromorphic function on $\bP^n$ is algebraic.  Let $P \colon \C \to \R$
be the defining polynomial of $S$.
Let $H = \overline{\{ z \mid P \circ F (z) = 0 \}}$.
We only need to show locally that $H$ is a subvariety and
that it is Levi-flat at all points of $H^*$.
Write $F$ in some set of affine coordinates as $F = f/g$
for two relatively prime polynomials $f$ and $g$.  If $P$ is a polynomial of
degree $d$ we notice that $\abs{g^d}^2 (P \circ F)$ is a polynomial
whose zero set
is precisely $H$ in the affine chart we have chosen.  Hence $H$
is a real-algebraic subvariety.  To see that it is Levi-flat, note
that locally it is always foliated by surfaces defined by the
set $\{ f = \lambda g \}$
for some constant $\lambda \in \C$.

To see the second part we refer to Example~\ref{ex:nodegen}.  
In that example we
construct an algebraic Levi-flat hypervariety such that there does not
exist a point contained in infinitely many leaves of the Levi
foliation.
If $H$ is defined by
a meromorphic function, there has to exist a point $p$ of indeterminacy since
the dimension is at least $2$.  
Define $f$ and $g$ in a given affine chart as above.  As 
the leaves of the Levi foliation are given by
$f(z) = \lambda g(z)$, we note that they all pass through $p$.
No such point $p$ exists on the hypervariety given in Example~\ref{ex:nodegen}.
\end{proof}

It is natural to ask the following question.  Can we define a Levi-flat
hypervariety by a meromorphic function as above, but choosing an arbitrary
real-analytic subset of $\C$ rather than an algebraic one.  The following
lemma says that this construction
would not yield a subvariety, or even a semianalytic set,
locally near a point
of indeterminacy of the meromorphic function.  
This lemma is the main new ingredient 
for the proof of Theorem~\ref{mainthmalt}.
If $H$ is defined by a meromorphic function as above,
then the function is
constant along the leaves of $H^*$.

\begin{lemma} \label{lfmeroalg}
Let $H \subset U \subset \C^n$, $n \geq 2$, be an irreducible
Levi-flat hypervariety
of $U$, and let
$p \in H$ be a point.
Suppose there exists a 
meromorphic function $F$ defined in $U$ such that $F$ is constant along
the leaves of $H^*$ and 
$p$ is a point
of indeterminacy of $F$.  Then there exists a one-dimensional algebraic subset
$S \subset \C$ such that $H \subset \overline{F^{-1}(S)}$.
\end{lemma}

\begin{proof}
First note that without loss of generality we can assume that $n=2$.
If we
pick a 2-dimensional subspace $V$ and find an $S$ such that
$\overline{F^{-1}(S)} \cap V$ contains $H \cap V$, then
since the inverse image of a single point under $F$ contains
the whole relevant leaf of the Levi foliation, 
$\overline{F^{-1}(S)}$ must then contain all of $H$ as $H$ is irreducible.

We can freely also pick a smaller neighborhood $U$ of $p$.  If the
conclusion of the lemma is true for a smaller neighborhood,
then it is true for the original $U$.
So
by perhaps picking a smaller $U$,
we can assume that the neighborhood $U$ is symmetric with respect to complex
conjugation and assume that $H$
complexifies to $U \times U$.
That is, the Taylor series of defining equation $\rho$ of $H$
converges on $U \times U$ if we replace $\bar{z}$ with a new variable $w$.

Let $F = f/g$ in $U$ where $f$ and $g$ are relatively prime.  If we look at the
map
\begin{equation}
\psi(z,w) := (f(z),g(z),\bar{f}(w),\bar{g}(w))
\end{equation}
and notice that it is a finite map
because $f^{-1}(0) \cap g^{-1}(0)$ must be a set of codimension 2, hence
a finite set.  If $\tilde{H} \subset
U \times U$ is the complexified $H$, then as $\psi$ is finite,
the image is also a complex subvariety.
We are really interested in $\varphi(H)$, where
\begin{equation}
\varphi(z) = (f(z),g(z)) .
\end{equation}
The image $\varphi(H)$ can be thought of as
a (possibly proper) subset of
$\psi(\tilde{H})$ intersected with the
totally real submanifold $\{ \zeta_1,\zeta_2,\zeta_3,\zeta_4 \mid
\zeta_1 = \bar{\zeta}_3,
\zeta_2 = \bar{\zeta}_4 \}$.  The point is that 
$\varphi(H)$
is semianalytic, that is,
near the origin
contained in a real-analytic subvariety $K$ of
the same dimension.

Notice that $\varphi(H)$ is Levi-flat
and $G(z) = z_1/z_2$ is constant along leaves of $\varphi(H)$.  That means that
$\varphi(H)$ contains complex lines through the origin.  
Take a defining function $r(z,\bar{z})$ for $K$.  Write
$r$ as
\begin{equation}
r(z,\bar{z}) = \sum_{j,k} r_{jk}(z,\bar{z})
\end{equation}
where $r_{jk}$ is homogeneous of order $j$ in $z$ and of order $k$
in $\bar{z}$.  Suppose that $z \in \varphi(H) \subset K$, then $\lambda z \in
\varphi(H) \subset K$
for some small open set of $\lambda$ and so
\begin{equation}
0 = \sum_{jk} r_{jk} (\lambda z,\bar{\lambda} \bar{z})
= \sum_{jk} \lambda^j \bar{\lambda}^k r_{jk} (z,\bar{z}) .
\end{equation}
By the same logic as in the proof of Proposition~\ref{basicfacts}
the set
$K$ is defined by a bihomogeneous polynomial.
In particular $K$ is real-algebraic.

We look at the algebraic set
\begin{equation}
\{ (z,\xi) \in \C^{2} \times \C 
\mid
z \in K, \text{ and }
\xi z_2 = z_1
\} .
\end{equation}
We project this set onto the $\xi$ variable.  By the theorem of
Tarski-Seidenberg, the projection must be semialgebraic.
It is not hard to see that the
set must be of dimension one.  A one-dimensional semialgebraic set
is contained in a one-dimensional algebraic set $S \subset \C$.  Hence,
$K \subset \overline{G^{-1}(S)}$, and as $\varphi(H) \subset K$
and as $F = G \circ \varphi$ then
$H \subset \overline{F^{-1}(S)}$.
\end{proof}

%%%%%%%%%%%%%%%%%%%%%%%%%%%%%%%%%%%%%%%%%%%%%%%%%%%%%%%%%%%%%%%%%%%%%%%%%%%%%

\section{Levi-flats and foliations} \label{lfandfolisec}

In this section we will prove Theorem~\ref{mainthmalt}.  If a Levi-flat
hypervariety of $\bP^n$ is locally defined by meromorphic functions and
has infinitely many compact leaves,
then it is algebraic and
furthermore defined by a global meromorphic function on $\bP^n$.  We will need 
the language of holomorphic foliations to prove this result.

A possibly singular holomorphic foliation $\sF$ of codimension one of a complex
manifold $M$ is given by
an open covering $\{ U_\iota \}$ with the following property.
In each $U_\iota$ there exists a holomorphic
one-form $\omega_\iota$ with $d \omega_\iota \wedge \omega_\iota = 0$.
If $U_\iota \cap U_\kappa \not= \emptyset$,
then $\omega_\iota$ and $\omega_\kappa$ must be proportional at every
point of the intersection.  A complex manifold is called a solution
if it satisfies the differential equation $\omega_\iota = 0$ in each
$U_\iota$.  The points where $\omega_\iota$ vanishes are called
the singular set of $\sF$ and denoted $\sing(\sF)$.  The set
$M \setminus \sing(\sF)$ is then a union of immersed complex hypersurfaces
called leaves of the foliation.
Note that the codimension
of the singularity of the foliation
can safely be taken to be at least 2, by dividing out
the coefficients of the form by a common divisor.
When talking about foliations of $\bP^n$,
we will say a leaf is compact
if its topological closure is of the same dimension.
In this case we will also use the word leaf
for the closure.
As we assume the singularity is of codimension at
least 2, a compact leaf is a complex analytic subvariety by the theorem
of Remmert and Stein.
See \cites{linsneto:note, CamachoNeto:book}
for more information on foliations in general.
All foliations in the sequel will be holomorphic of codimension one.

The Levi foliation of a Levi-flat hypervariety
does not necessarily extend (even locally) to a foliation of a neighborhood
of the hypervariety, at least not in the above sense,
see Brunella~\cite{Brunella:lf}.
If the Levi-flat hypervariety is such that
locally there exists a meromorphic function
$F=f/g$ (in lowest terms) that is constant along leaves of $H^*$,
then the foliation extends locally.
The leaves
are defined by components of the sets $\{ f=\lambda g \}$ for a constant
$\lambda$ and
the form is given by $\omega = f(dg) - g(df)$.

Of course, the condition that the foliation extends is a
necessary condition for a hypervariety to be defined in the same manner as
in Theorem~\ref{mainthmalt}.
If we further know that the leaves of the Levi foliation are compact, these
two conditions
turn out to be sufficient.

As we said in the introduction, we will prove
the theorem for semianalytic sets.  We also need not require that
the foliation extending that of $H$ be locally first integrable.
The main feature of semianalytic sets we will use is that at
each point, a germ of a semianalytic set is contained in a germ of
a real-analytic set of the same dimension.

\begin{thm} \label{samainthm}
Let $H \subset \bP^n$, $n \geq 2$, be a connected semianalytic set of
real dimension $2n-1$.  Suppose that $H = \bigcup_\iota L_\iota$,
where $L_\iota$ are complex analytic
hypervarieties of $\bP^n$.
Assume that for each $p \in H$, there
exists a neighborhood $U$ of $p$ and a holomorphic
foliation on $U$ such that $L_\iota \cap U$ are invariant.

Then, there exists a global rational function $R \colon \bP^n \to \C$
and
a real-algebraic one-dimensional
subset $S \subset \C$ such that $H \subset
\overline{R^{-1}(S)}$.
\end{thm}

Essentially we are asking for a foliation of a neighborhood of $H$, and $H$
to be an invariant set of the foliation.
To extend the foliation into all of $\bP^n$ we use the following result of
Lins Neto~\cite{linsneto:note}.

\begin{thm}[Lins Neto]  \label{folext}
Let $M$ be a Stein manifold, $\dim(M) \geq 2$.  Let $K \subset M$ be
compact, with $M \setminus K$ connected,
and let $\sF$ be a singular holomorphic foliation of $M \setminus K$
where $\codim(\sing(\sF)) \geq 2$.  Then $\sF$ extends to a singular
holomorphic foliation on $M$.
\end{thm}

To be able to apply Theorem~\ref{folext} we need to find Stein manifolds
inside $\bP^n$.

\begin{thm}[Takeuchi~\cite{Takeuchi}] \label{Takeuchithm}%, Kiselman?
Let $U \subset \bP^n$ be an open set such that $U \not= \bP^n$.  Suppose
that $U$ is pseudoconvex (satisfies Kontinuit\"atssatz), then $U$ is Stein.
\end{thm}

%As $H$ has infinitely many compact leaves, at least one of
%them, call it $L$, lies inside $\overline{H^*}$.

Take one complex variety $L = L_\iota$ that lies in $H$.
The set $\bP^n \setminus L$
is Stein.  If we have a foliation of a neighborhood of $H$,
we have a foliation of a neighborhood of $L$.
We can then apply Theorem~\ref{folext} to
get a foliation of $\bP^n$.
The leaves of the foliation must coincide with
the complex varieties near $L$ that are part of $H$.
Note that $H$ is connected, so the $L_\iota$ are leaves of
the extended foliation of $\bP^n$.

Once we have the foliation extended to all of $\bP^n$, we will need to
find the rational function $R$.  Therefore,
we apply the following classical theorem
of Darboux (see \cite{Ince} page 29)
generalized by Jouanolou, see~\cite{jou:pfaff} Theorem 3.3 page 102.

\begin{thm}[Darboux-Jouanolou] \label{darboux}
If $\sF$ is a singular holomorphic foliation on $\bP^n$ with infinitely
many compact leaves, then $\sF$ has a rational first integral.
\end{thm}

As $H$ is of dimension $2n-1$, it must contain infinitely many
complex varieties.  As these coincide with the leaves of $\sF$ (the extended
foliation), $\sF$ has infinitely many compact leaves and hence has a rational
first integral.
Therefore, the final piece
of the proof of Theorem~\ref{mainthmalt} is the following
lemma, which may be of independent
interest.

\begin{lemma} \label{finalpiece}
Suppose that there exists a
singular holomorphic foliation $\sF$ of $\bP^n$, $n \geq 2$,
with a rational first integral $R$.
Let $H \subset \bP^n$ be a
connected semianalytic set of real dimension $2n-1$
that is an invariant set of $\sF$.
Then there exists 
a real-algebraic one-dimensional
subset $S \subset \C$ such that $H \subset
\overline{R^{-1}(S)}$.
\end{lemma}

\begin{proof}
We can assume that $H$ is closed just by taking the closure, which is also
semianalytic and invariant.
We can write
$H = \bigcup_\iota L_\iota$, where $L_\iota$ are irreducible
complex analytic hypervarieties (leaves of $\sF$).
Let $R$ be the first integral of $\sF$.  $R$ must be constant 
along the $L_\iota$.
We find a point $p \in \bP^n$ that is
a point of indeterminacy for $R$, which exists because $n \geq 2$.
Further a point $p$ of indeterminacy has to lie on $H$,
since it must be in the closure of
every leaf $L_\iota$.  That is, write $R=f/g$, then without loss of
generality there is some fixed $\lambda \not= 0$, such that
$f/g = \lambda$ on $L_\iota$.  The numerator $f$ must be
zero on $L_\iota$ and then $g$ must also be zero at the same point and that
must be a point of indeterminacy.

We can apply Lemma~\ref{lfmeroalg} near $p$ to find the required $S$.  As
$H \subset \overline{R^{-1}(S)}$ locally near $p$, $H$ is a union of the
$L_\iota$, and $p \in L_\iota$ for all $\iota$, then 
$H \subset \overline{R^{-1}(S)}$.
\end{proof}

\begin{proof}[Proof of Theorem~\ref{samainthm}]
For every point $p \in H$ we have a neighborhood $U$ and a 
foliation on $U$ extending the Levi foliation of $H^* \cap U$.
We call $H^*$ the smooth part of real dimension $2n-1$
just like for hypervarieties.

Suppose we have two connected neighborhoods $U_1$
and $U_2$ such that $U_1 \cap U_2$ is nonempty and connected.
Further, assume there exist holomorphic one-forms
$\omega_j$ on $j=1,2$ that define a foliation extending the foliation of
$H^*$.  If we can show that $\omega_1$ is proportional to $\omega_2$
then we have a foliation of $U_1 \cap U_2$.  Take a small neighborhood
of some point $p \in H^*$.  We know that $\omega_1$ must be proportional to
$\omega_2$ for all points of $H^*$ near $p$.  $H^*$ is a real hypersurface,
thus they
are proportional in a whole neighborhood as they are holomorphic.  Since
$U_1 \cap U_2$ is connected, we are done by analytic continuation.

We can choose a covering of $H$ that satisfies the
above conditions
for every pair of intersecting neighborhoods.
Hence, if the foliation of $H^*$ extends locally near every point of
$H$, then there exists a neighborhood $U$ of $H$ and
a foliation $\sF$ on $U$ that extends the foliation of $H^*$.
Again, we can assume that the codimension of the singularity
of the foliation is at least two.

We pick one complex hypervariety $L$ that lies in $H$
and we apply Theorem~\ref{folext} to
extend
the foliation to a foliation on all of $\bP^n$.
As we said above, the foliation $\sF$ has infinitely many compact leaves.
We can apply Theorem~\ref{darboux} to get a rational first integral.

Finally we appeal to Lemma~\ref{finalpiece} which has the same conclusion
as our theorem.
\end{proof}

Once we have Theorem~\ref{samainthm}, it is not too hard to finish
the proof of 
Theorem~\ref{mainthmalt}.  We can notice that 
$\overline{H^*}$ is semianalytic, and we just need to show that it is
a union of compact leaves.  However, it is easier to modify
the above proof.  As above, we need not require the foliation 
to be locally first integrable.

\begin{thm} \label{mainthmalt2}
Let $H \subset \bP^n$, $n \geq 2$, be an irreducible
Levi-flat hypervariety with infinitely many compact leaves.
Assume that for each $p \in \overline{H^*}$, there
exists a neighborhood $U$ of $p$ and a holomorphic
foliation on $U$ extending the Levi foliation of $H^*$.

Then, there exists a global rational function $R \colon \bP^n \to \C$
and
a real-algebraic one-dimensional
subset $S \subset \C$ such that $H \subset
\overline{R^{-1}(S)}$.
In particular, $H$ is semialgebraic; it is contained in an
algebraic Levi-flat
hypervariety.
\end{thm}

\begin{proof}%[Proof of Theorem~\ref{mainthmalt}]
Follow the same argument as in the proof of Theorem~\ref{samainthm}.
By exactly the same argument, we have a foliation of a
neighborhood of all of $\overline{H^*}$.
If there are infinitely many
compact leaves of $H$, one is contained in $\overline{H^*}$.  We have
a foliation of a neighborhood of this compact leaf and we can extend the
foliation to a foliation of $\bP^n$.  So
we have a foliation
of $\bP^n$ that extends the foliation of $H^*$.  As it has infinitely
many compact leaves, it has a first integral $R$.
The set $\overline{H^*}$ is invariant and
satisfies the hypothesis of
Lemma~\ref{finalpiece}.
\end{proof}

%%%%%%%%%%%%%%%%%%%%%%%%%%%%%%%%%%%%%%%%%%%%%%%%%%%%%%%%%%%%%%%%%%%%%%%%%%%%%

\section{Extending foliations}
\label{extfolsec}

We have already proved Theorem~\ref{mainthmalt}, but
it will be interesting to also
prove the following stronger result about foliations, which does not use
the compactness of leaves.  This result
is also of independent interest and is
essentially an extension of a similar result by Lins-Neto to singular
Levi-flat hypervarieties.

\begin{thm} \label{mainfolthm}
Suppose $H \subset \bP^n$, $n \geq 2$, is an irreducible
Levi-flat hypervariety.
Assume that for each $p \in \overline{H^*}$, there
exists a neighborhood $U$ of $p$ and a meromorphic function $F$ defined
on $U$ such that $F$ is constant along leaves of $H^*$.

Then, there exists a singular holomorphic foliation $\sF$ of $\bP^n$
that agrees with the foliation of $H^*$.
\end{thm}

We already know we have a foliation of a neighborhood of $\overline{H^*}$.
We notice the following corollary of the theorem of Takeuchi.  Once the
following lemma is proved, the proof of Theorem~\ref{mainfolthm} follows at
once. 

\begin{lemma} \label{compstein}
Let $H \subset \bP^n$ be a Levi-flat
hypervariety, such that for every $p \in \overline{H^*}$ there exists
a neighborhood $U$ and a meromorphic function $F$ such that $F$ is
constant along leaves of $H^*$.
Then all the connected components of $\bP^n \setminus \overline{H^*}$
are Stein.
\end{lemma}

This corollary follows after we have shown that through every point of
$\overline{H^*}$ there exists a germ of a complex hypervariety contained in
$\overline{H^*}$, hence $\overline{H^*}$ is pseudoconvex from all sides.
A weaker theorem, that through every point of
$\overline{H^*}$ there exists a complex hypervariety contained in $H$
was essentially proved by Fornaess (see~\cite{kohn:subell} Theorem 6.23).
The statement by Fornaess assumes that $H$ is nonsingular, but that is
not used in the proof.
See also
Burns and Gong~\cite{burnsgong:flat} for more information regarding this
point.

If $H$ was not a complex hypersurface near any point, we would be done by the
theorem of Fornaess.  However,
Example~\ref{cartanumbex} shows that it is possible to have
an irreducible Levi-flat hypervariety with a component
that is a complex hypervariety.  We will prove the following lemma,
which, together with Takeuchi's theorem,
implies Lemma~\ref{compstein}, and hence Theorem~\ref{mainfolthm}.

\begin{lemma}
Let $H \subset U \subset \C^k$ be a 
Levi-flat hypervariety and $p \in
\overline{H^*}$.  Suppose there exists a meromorphic function $F$
defined in $U$ that is constant along leaves of $H^*$.  Then
there exists a germ of a complex hypervariety $(L,p)$ such that
$(L,p) \subset (\overline{H^*},p)$.
\end{lemma}

\begin{proof}
First let us assume that $p$ is a point of indeterminacy of $F$.
Let $F = f/g$ written in lowest terms.
We can follow the proof of Lemma~\ref{lfmeroalg} to note that
the image of $H$ under the map $z \mapsto (f(z),g(z))$ is a complex
cone.  Therefore, given a constant $\lambda$, set
$\{ f = \lambda g \}$ contains a leaf
of $H$ going through the origin.  I.e.\@ there are infinitely many leaves
of $H$ going through the origin.  Only finitely many leaves
can form a ``stick''
of an umbrella, and hence infinitely many are contained in $\overline{H^*}$.

If $p$ is not a point of indeterminacy of $F$, we can assume $F$ is
holomorphic.  By taking $U$ smaller, we could assume $F$ is holomorphic in all
of $U$ and assume $F(p)=0$.  Define the graph
$\Gamma_F := \{ (z,\xi) \mid \xi = F(z) \}$.  After a possible linear change of
coordinates in the $z$ variable
we can apply the Weierstrass preparation theorem to get $\Gamma_F$
defined by
\begin{equation}
z_k^d + \sum_{j=0}^{d-1} a_j(z',\xi) z_k^j = 0 ,
\end{equation}
where $z' = (z_1,\ldots,z_{k-1})$.  The set
$V_\lambda := \{ z \mid F(z) = \lambda \}$ is a multigraph over the $z'$
of multiplicity at most $d$.  That is, we have a holomorphic
function $h \colon U' \subset \C^{k-1} \to \C_{sym}^d$ (a multifunction)
and $V_\lambda$ is the set $\{ z \mid z_k \in h(z') \}$.  See
\cite{Whitney:book} for more information on symmetric powers and 
multifunctions.

Pick a sequence of $\lambda_j \to 0$, such that $V_{\lambda_j}$ contains
a branch $V'_{\lambda_j} \subset H$.
As there can locally be at most finitely many branches of $H$
that are complex hypervarieties, we can assume that
$V'_{\lambda_j} \subset \overline{H^*}$ for all $j$.  As $V_\lambda$ is a multigraph of
multiplicity $d$, there must exist a single integer $m$ such that
each $V'_{\lambda_j}$ is a multigraph of multiplicity $m$.
Assume $V'_{\lambda_j}$ is the multigraph of
$h_j \colon U' \subset \C^{k-1} \to \C_{sym}^{m}$.
The functions $h_j$ are
bounded and hence we can pass to a convergent subsequence.
That is, there exists a complex hypervariety $V$ that is the limit of
$V'_{\lambda_j}$.  Since
$V'_{\lambda_j} \subset \overline{H^*}$ then $V \subset \overline{H^*}$,
furthermore, $p \in V$ as $\lambda_j \to 0$ and $F(p) = 0$.
\end{proof}

%It is a natural question and perhaps not an one to answer, weather a
%foliation with an invariant real-analytic subvariety or
%an invariant semianalytic set of real dimension $2n-1$
%necessarily has compact leaves.  It seems reasonable to
%conjecture that the answer is negative.

%%%%%%%%%%%%%%%%%%%%%%%%%%%%%%%%%%%%%%%%%%%%%%%%%%%%%%%%%%%%%%%%%%%%%%%%%%%%%

\section{Nonalgebraic hypervarieties with compact leaves}
\label{compactleavessec}

In this paper,
we have mostly studied Levi-flat hypervarieties (or semianalytic sets)
with compact leaves.  Each compact leaf is algebraic.  Therefore, the
following construction gives the most obvious type of Levi-flat hypersurface
with compact leaves.  For $z \in \C^{n+1}$ let
\begin{equation}
f(z,t) = \sum_{\abs{\alpha} = d} c_{\alpha}(t) z^{\alpha} ,
\end{equation}
where $c_\alpha(t)$ are real-analytic functions of 
$t \in (a,b) \subset \R$.  If $c_\alpha$ are analytic up to $a$ and $b$,
then
\begin{equation} \label{clform1}
H = \{ z \in \C^{n+1} \mid f(z,t) = 0, \text{ for some $t \in (a,b)$} \} 
\end{equation}
is a subanalytic Levi-flat hypersurface, which is a complex cone.  Hence $\sigma(H)$ is a
subanalytic Levi-flat hypersurface in $\bP^n$.
%It is generally not easy to
%check if $\sigma(H)$ is semianalytic or perhaps even contained in a
%real subvariety of $\bP^n$.

Define the function of $(w,t) \in \C^{N+1} \times (a,b)$ by
\begin{equation}
F(w,t) = \sum_{k=1}^{N+1} c_{\alpha}(t) w_k .
\end{equation}
The set
\begin{equation} \label{clform2}
H' = \{ w \in \C^{N+1} \mid F(z,t) = 0, \text{ for some $t \in (a,b)$} \} 
\end{equation}
is a subanalytic Levi-flat hypersurface whose leaves are complex
hyperplanes.  As before $\sigma(H') \subset \bP^n$ is also subanalytic Levi-flat
hypersurface.
Let $\sZ \colon \C^{n+1} \to \C^{N+1}$, where $N+1$
is the number of distinct degree $d$ monomials, be the degree $d$
Veronese mapping.
That is,
$\sZ$ is
the mapping $z \mapsto \bigoplus_{\abs{\alpha}=d} z^\alpha$.
We then have
\begin{equation}
H = \sZ^{-1} (H') .
\end{equation}
Therefore, to study hypersurfaces of the form \eqref{clform1} we need only
study Levi-flat hypersurfaces of the form \eqref{clform2} with leaves
being complex hyperplanes.

\begin{example} \label{ex:nonalg}
Let us build a semianalytic Levi-flat hypersurface of $\bP^2$
with compact leaves, which is a
small perturbation of an algebraic Levi-flat hypervariety of $\bP^2$, but is
not algebraic itself.  This example suggests that any analogue of Chow's
theorem for Levi-flat hypervarieties will likely have to require compact
leaves.

First let us construct the algebraic Levi-flat hypervariety.  Take 
\begin{equation}
H = 
\{ z \in \C^3 \mid z_1 + x z_2 + y z_3 = 0, x^2+y^2=1, x \in \R, y \in \R \}.
\end{equation}
That is, $H$ is the projection of a variety in $\C^3 \times \R^2$ onto
$\C^3$.  It is not just semialgebraic, it is in fact a real
hypervariety in $\C^3$, and of course Levi-flat with leaves that are complex 
hyperplanes.  To see that $H$ is a variety, write
$z_j = s_j + i t_j$.  We have 
$s_1 + x s_2 + y s_3 = 0$ and
$t_1 + x t_2 + y t_3 = 0$.  Solve for $x$ and $y$ to get
\begin{equation}
x=-\frac{s_3 t_1-s_1 t_3}{s_3 t_2-s_2 t_3},
\quad
y=\frac{s_2 t_1-s_1 t_2}{s_3 t_2-s_2 t_3}.
\end{equation}
Therefore, $H$ is defined by
\begin{equation}
(s_3 t_1-s_1 t_3)^2
+
(s_2 t_1-s_1 t_2)^2
= (s_3 t_2-s_2 t_3)^2 .
\end{equation}

This equation defines a Levi-flat complex cone in $\C^3$ and hence a Levi-flat
hypervariety in $\bP^2$.

To define a perturbation of $H$, we want to perturb $x^2+y^2=1$.  Suppose
we take a real-analytic $f(x)$ that is a small perturbation of
$x$, and such that
$C = \{ \R^2 \mid f(x)^2+y^2=1 \}$
is not contained in an algebraic curve.
We can also ensure
that near each point $p$ on the curve $C$
we can parametrize $C$ by a one-to-one real-analytic $\gamma \colon
(-\epsilon,\epsilon) \to C$, and we only need to pick finitely many such
$\gamma$'s to parametrize all of $C$.
That is, $C$ is a compact topological manifold.

We now need only show that 
\begin{equation}
H' = 
\{ z \in \C^3 \mid z_1 + x z_2 + y z_3 = 0, f(x)^2+y^2=1, x \in \R, y \in \R \}
\end{equation}
is semianalytic for all $p \in H'$ except $p=0$.  Then we need to show
that $H'$ is not contained in a real-algebraic variety.
We then obtain a Levi-flat semianalytic set in
$\bP^2$ with compact leaves that is not contained in a real-algebraic
Levi-flat.

Define $H'' \subset \C^3 \times \R^2$ by
$z_1 + x z_2 + y z_3 = 0$ and $f(x)^2+y^2=1$.
Take a point $(\xi_1, \xi_2, \xi_3, x_0, y_0) \in H''$ such
that $(\xi_1,\xi_2,\xi_3) \not= (0,0,0)$.  Find a $\gamma =
(\gamma_1,\gamma_2)$ as above
for the $C$ near the point $(x_0,y_0)$.
The function
\begin{equation} \label{parameq}
t \mapsto
\xi_1 + \gamma_1(t) \xi_2 + \gamma_2(t) \xi_3
\end{equation}
is not identically zero.  Hence, we can apply
Weierstrass preparation theorem to
the function
$z_1 + \gamma_1(t) z_2 + \gamma_2(t) z_3$
of $(z,t)$ with respect to the $t$ variable.
The projection of $H'' \cap U$ to $\C^3$ for
some neighborhood $U$ of $(\xi_1, \xi_2, \xi_3, x_0, y_0)$
is the same as the projection of
$\{ z_1 + \gamma_1(t) z_2 + \gamma_2(t) z_3 = 0 \}$
for some small interval of $t$.  If we know that the projection
of this set to $\C^3$ is semianalytic, we are done.

The above claim is achieved by the
following version of Tarski-Seidenberg theorem by \L{}ojasiewicz,
see Theorem 2.2 in~\cite{BM:semisub}.
Let us set up some terminology.
Suppose $\sA(U)$ is any ring of real valued functions on an open set
$U \subset \R^n$.
Define $\sS(\sA(U))$ to be the smallest
set of subsets of $U$, which contain the sets
$\{ x\in U \mid f(x) > 0 \}$ for all $f \in \sA(U)$,
and is closed under finite union, finite intersection and complement.
A set $V \subset \R^n$ is semianalytic
if and only if for each $x \in \R^n$, there exists a neighborhood
$U$ of $x$, such that $V \cap U \in \sS(\sO(U))$, where $\sO(U)$ 
denotes the real-analytic real valued functions.
Let $\sA(U)[t]$ denote the ring of polynomials in $t \in \R^m$
with coefficients in $\sA(U)$.

\begin{thm}[Tarski-Seidenberg-{\L}ojasiewicz] \label{thm:tsloj}
Suppose that $V \subset U \times \R^m \subset \R^{n+m}$,
is such that $V \in \sS(\sA(U)[t])$.
Then the projection of $V$ onto the first $n$ variables
is in $\sS(\sA(U))$.
\end{thm}

Consequently, if we can locally Weierstrass the defining function with
respect to the $t$ variable, we can project onto the remaining variables and
obtain a semianalytic set.  Of course, the Weierstrass theorem will only
apply in some neighborhood, and hence for a small finite interval of the
$t$.  We only need to do the projection for $t$ in a compact
interval for finitely many curves $\gamma$.  A finite union of semianalytic
sets is semianalytic.

Finally we must show that $H'$ is not contained in a real-algebraic
hypervariety.  Fix $z_2 = -1$ and $z_3 = -i$.  The defining equations become
\begin{equation}
z_1 = x + iy, \quad f(x)^2 +y^2 = 1 .
\end{equation}
We picked $f(x)$ precisely in such a way that this set projected onto
$z_1$ is not contained
in an algebraic curve.
\end{example}

%%%%%%%%%%%%%%%%%%%%%%%%%%%%%%%%%%%%%%%%%%%%%%%%%%%%%%%%%%%%%%%%%%%%%%%%%%%%%

%FIXME: else I don't get links, weird
%\renewcommand\MR[1]{\relax\ifhmode\unskip\spacefactor3000 \space\fi
  %\def\@tempa##1:##2:##3\@nil{%
    %\ifx @##2\@empty##1\else\textbf{##1:}##2\fi}%
  %\href{http://www.ams.org/mathscinet-getitem?mr=#1}{MR \@tempa#1:@:\@nil}}
\def\MR#1{\relax\ifhmode\unskip\spacefactor3000 \space\fi%
  \href{http://www.ams.org/mathscinet-getitem?mr=#1}{MR#1}}

\end{document}